\documentclass[a4, 12pt]{amsart}
\usepackage[mathscr]{eucal}
\usepackage{amssymb}
\usepackage{latexsym}
\usepackage{amsthm}
\usepackage{color}
\theoremstyle{plain}
\newtheorem{theorem}{Theorem}[section]

\newtheorem{remark}{Remark}[section]
\newtheorem{lemma}{Lemma}[section]
\newtheorem{proposition}{Proposition}[section]

\makeatletter
\@addtoreset{equation}{section}

\title[Rigidity theorems of complete self-shrinkers in the Euclidean space]
{Rigidity theorems of complete self-shrinkers in the Euclidean space}

\author [Z. Li, G. Wei and C. Xing]{Zhi Li, Guoxin Wei and Cheng Xing}

\address{Zhi Li \\  \newline \indent School of Mathematics and Statistics, Henan Normal University,
\newline \indent 453007, Xinxiang, Henan, China.}
\email{lizhihnsd@126.com}

\address{Guoxin Wei \\ \newline \indent School of Mathematical Sciences, South China Normal University,
\newline \indent 510631, Guangzhou,  China.}
\email{weiguoxin@tsinghua.org.cn}

\address{Cheng Xing\\  \newline \indent School of Mathematics and Statistics, Henan Normal University,
\newline \indent 453007, Xinxiang, Henan, China.}
\email{xingcheng@htu.edu.cn}

\begin{document}
\maketitle

\begin{abstract}
We classify complete Lagrangian self-shrinking surfaces in $\mathbb{R}^{4}$ with constant squared norm of the mean curvature vector, without assuming an upper bound for the second fundamental form. Their images are a plane through the origin, the round cylinder $S^1(1)\times \mathbb{R}$, or the Clifford torus $S^1(1)\times S^1(1)$. More generally, a complete self-shrinking surface in $\mathbb{R}^{4}$ with constant $|\vec{H}|^2$ is either a plane, a round cylinder, or a complete minimal immersion into $S^3(\sqrt2)$. The proof uses the function $|x|^2-|\vec{H}|^2$, the Bochner formula, and completeness.
\end{abstract}


\section{introduction}
\vskip2mm
\noindent

An immersion $x:M^n\to\mathbb{R}^{n+p}$ is called a self-shrinker if
\begin{align}\label{1.1-1}
\vec{H}+x^{\perp}=0,
\end{align}
where $\vec{H}=trace B$ is the mean curvature vector and $x^{\perp}$ is the normal component of the position vector. We use the trace convention for mean curvature, so that $\Delta x=\vec{H}$. Unless otherwise stated, all manifolds are smooth, connected, and without boundary, and completeness refers to the induced Riemannian metric.

Self-shrinkers generate homothetically shrinking solutions of mean curvature flow and arise in the analysis of Type I singularities \cite{Hui90}. They are stationary for the Gaussian weighted area $$\int_{M}e^{-\frac{|x|^{2}}{2}}d\sigma$$
under compactly supported variations. Equivalently, they are minimal submanifolds of Euclidean space equipped with the conformal metric $\bar g=e^{-|x|^2/n}g_{\mathrm{Euc}}$. Classification and rigidity results have been studied extensively in both codimension one and higher codimension; see \cite{CL13, CP15, CW15, CZ13, CW, DX13, DX14, LN11, LW14, Smo05, XX17} and the references therein.

The Lagrangian condition is preserved under smooth Lagrangian mean curvature flow. This motivates the study of rigidity and classification for Lagrangian self-shrinkers. Castro and Lerma \cite{Cas} classified Hamiltonian stationary Lagrangian self-shrinkers in $\mathbb{R}^4$ and subsequently established rigidity results for the Clifford torus among compact Lagrangian self-shrinkers \cite{Cas1}. Li and Wang \cite{LW3} obtained further characterizations of the Clifford torus. Cheng, Hori, and Wei \cite{CHW} classified complete Lagrangian self-shrinkers in $\mathbb{R}^4$ with constant squared norm of the second fundamental form. For related constructions and developments, see \cite{Anc06, LW10, Neves11, Smo00}; related rigidity results for Lagrangian submanifolds and $\xi$-submanifolds appear in \cite{LW2,LC16}.

In particular, Castro and Lerma \cite{Cas1} characterized the Clifford torus among compact Lagrangian self-shrinkers with constant $|\vec{H}|^2$. Li, Wang and Wei \cite{LWW} established the corresponding classification in the complete case under an upper bound for $|B|^2$, using a generalized maximum principle. The purpose of this paper is to remove that curvature bound and to give a direct proof. Our first result is the following.

\begin{theorem}\label{theorem 1.1}
Let $x:M^{2}\to\mathbb{C}^{2}=\mathbb{R}^{4}$ be a complete connected Lagrangian self-shrinker with constant $|\vec{H}|^{2}$. Then $x(M^{2})$ is one of
$$
\mathbb{R}^{2},\quad S^{1}(1)\times\mathbb{R},\quad S^{1}(1)\times S^{1}(1).
$$
The plane passes through the origin, the axis of the cylinder passes through the origin, and the torus is centered at the origin.
\end{theorem}

The same argument also yields a classification without the Lagrangian assumption.

\begin{theorem}\label{theorem 1.2}
Let $x:M^{2}\to\mathbb{R}^{4}$ be a complete connected self-shrinker with constant $c=|\vec{H}|^2$. Then $c\in\{0,1,2\}$, and exactly one of the following occurs:
\begin{enumerate}
\item $c=0$: $x(M^{2})$ is a plane through the origin;
\item $c=1$: $x(M^{2})$ is the round cylinder $S^{1}(1)\times\mathbb{R}$ in a linear $\mathbb{R}^3\subset\mathbb{R}^{4}$;
\item $c=2$: $x(M^{2})$ is a minimal surface into $S^{3}(\sqrt{2})$.
\end{enumerate}
\end{theorem}

\begin{remark}
In fact, for codimension greater than $2$ in Theorem \ref{theorem 1.2}, the conclusion still holds.
\end{remark}

\vskip5mm
\section {Preliminaries}
\vskip2mm

\noindent

Let $x:M^{n}\to \mathbb{R}^{n+p}$ be an isometric immersion. Denote by $D$, $\nabla$ and $\nabla^{\perp}$ the Euclidean, Levi-Civita, and normal connections, respectively. Write $x=T+x^{\perp}$, where $T=x^{\top}$. For tangent vector fields $X,Y$ and a normal vector field $\eta$, the Gauss and Weingarten formulas are
$$D_{X}Y=\nabla_{X}Y+B(X,Y), \quad D_{X}\eta=-A_{\eta}X+\nabla^{\perp}_{X}\eta,$$
where $\langle A_{\eta}X, Y\rangle=\langle B(X,Y), \eta\rangle$.
Taking tangential and normal components of $D_{X}x=X$ gives
\begin{equation}\label{2.1-1}
\nabla_{X} T=X+A_{x^\perp} X, \quad \nabla^{\perp}_{X}x^{\perp}=-B(X, T).
\end{equation}
Define the drift Laplacian by
\begin{equation*}
\mathcal{L}f=\Delta f-\langle T,\nabla f\rangle.
\end{equation*}
where $\Delta$ and $\nabla$ denote the Laplacian and the gradient
operator, respectively.

\noindent
When \eqref{1.1-1} holds, \eqref{2.1-1} becomes
\begin{equation*}
\nabla_{X} T=X+A_{\vec{H}} X, \quad \nabla^{\perp}_{X}\vec{H}=B(X, T).
\end{equation*}

At an arbitrary point $p$, choose a geodesic orthonormal frame, i.e., $\nabla_{e_{i}} e_{j}(p) = 0$. Denote
$$
B_{ij} = B(e_{i}, e_{j}), \qquad T_{j} = \langle T, e_{j} \rangle.
$$
Then
$$
\nabla^{\perp}_{e_{i}}\vec{H}=\sum_{j}B_{ij} T_{j}, \quad \langle\nabla_{e_{i}} T, e_{j}\rangle = \delta_{ij}-\langle B_{ij}, \vec{H} \rangle.
$$
Differentiation yields
\begin{equation}\label{2.1-2}
\nabla^{\perp}_{e_{j}}\nabla^{\perp}_{e_{i}}\vec{H}=\sum_{k}(\nabla^{\perp}_{j} B)_{ik} T_{k}+B_{ij}-\langle B_{jk}, \vec{H} \rangle B_{ik}.
\end{equation}
Tracing and applying the contracted Codazzi equation gives
\begin{equation*}
\begin{aligned}
\Delta^{\perp} \vec{H}&= \sum_{i,j} (\nabla^{\perp}_{i} B)_{ij} T_{j} + \sum_{i,j} B_{ij} \big(\delta_{ij}-\langle B_{ij}, \vec{H} \rangle \big) \\
&=\nabla^{\perp}_{T}\vec{H}+\vec{H}-\sum_{i,j}\langle B_{ij}, \vec{H}\rangle B_{ij}.
\end{aligned}
\end{equation*}
Since $\sum_{i,j} \langle B_{ij}, \vec{H}\rangle^{2}=|A_{\vec{H}}|^{2}$, the product rule implies
\begin{equation}\label{2.1-3}
\frac{1}{2}\mathcal{L} |\vec{H}|^{2}=|\nabla^{\perp}\vec{H}|^{2}+|\vec{H}|^{2}-|A_{\vec{H}}|^{2}.
\end{equation}

For completeness, we also record the Simons identity for $B$ (cf.\ \cite{DX08}). Let $\{e_{\alpha}\}$ be a local orthonormal normal frame, and define
$$S_{\alpha\beta}=\sum_{i,j}\langle B_{ij}, e_{\alpha}\rangle\langle B_{ij}, e_{\beta}\rangle.$$
Then
\begin{equation}\label{2.1-4}
\begin{aligned}
\frac{1}{2}\Delta|B|^{2}=&|\nabla^{\perp} B|^{2}+\langle\nabla^{\perp}_{e_{j}}\nabla^{\perp}_{e_{i}}\vec{H}, B_{ij}\rangle+\langle\vec{H}, B_{ij}\rangle\langle B_{ik}, B_{jk}\rangle\\
&-\sum_{\alpha\neq\beta}|[A_{e_{\alpha}}, A_{e_{\beta}}]|^{2}-\sum_{\alpha, \beta}S^{2}_{\alpha\beta}.
\end{aligned}
\end{equation}
Substitution of \eqref{2.1-2} into \eqref{2.1-4}, together with Codazzi symmetry, gives
\begin{equation}\label{2.1-5}
\frac{1}{2}\mathcal{L} |B|^{2}=|\nabla^{\perp} B|^{2}+|B|^{2}-\sum_{\alpha\neq\beta}|[A_{e_{\alpha}}, A_{e_{\beta}}]|^{2}-\sum_{\alpha, \beta}S^{2}_{\alpha\beta}.
\end{equation}
The contractions can equivalently be written as
$$\sum_{\alpha\neq\beta}|[A_{e_{\alpha}}, A_{e_{\beta}}]|^{2}=2\sum_{i,j,k,l}\Big(\langle B_{ij}, B_{il}\rangle\langle B_{jk}, B_{kl}\rangle-2\langle B_{ij}, B_{kl}\rangle\langle B_{il}, B_{jk}\rangle\Big),$$ $$\sum_{\alpha, \beta}S^{2}_{\alpha\beta}=\sum_{i,j,k,l}\langle B_{ij}, B_{kl}\rangle^{2}.$$
The proof below only requires \eqref{2.1-1}--\eqref{2.1-5}. For surfaces, we will also use the Gauss equation
$$2K=|\vec{H}|^2-|B|^2.$$

\begin{lemma}\label{lemma 2.1}
Let $(M, g)$ be a nonempty complete Riemannian manifold, and let $u\geq0$ be a smooth function satisfying $|\nabla u|^{2}=4u$. Then $u^{-1}(0)\ne\varnothing$. For each $p\in\Omega=\{u>0\}$, there is a unit-speed curve starting at $p$, of length $\sqrt{u(p)}$, that converges to a point of $u^{-1}(0)$.
\end{lemma}
\begin{proof}
If $u\equiv0$, there is nothing to prove. On $\Omega$, set $r=\sqrt{u}$. Then $|\nabla r|=1$. For $p\in\Omega$, let $\gamma:[0, \tau)\to\Omega$ be the maximal forward integral curve of $-\nabla r$ with $\gamma(0)=p$. Along this curve,
$$
|\dot{\gamma}|=1,\quad \frac{d}{dt}r(\gamma(t))=-1,\quad
r(\gamma(t))=r(p)-t.
$$
Thus $\tau\leq r(p)<\infty$. Moreover,
$d(\gamma(s),\gamma(t))\leq|s-t|$, so completeness implies that $\gamma(t)$ converges to a point $q\in M$ as $t\uparrow\tau$.
If $\tau<r(p)$, continuity gives $u(q)=(r(p)-\tau)^2>0$, so $q\in\Omega$. The vector field $-\nabla r$ is smooth near $q$, and local existence for ordinary differential equations extends $\gamma$ beyond $\tau$, contradicting maximality. Hence $\tau=r(p)$ and $u(q)=0$. The length of $\gamma$ is exactly $r(p)$.
\end{proof}

 \vskip10mm
\section{Proofs of the rigidity theorems}

\vskip2mm
\begin{proposition}\label{proposition 3.1}
Let $x:M^{2} \to \mathbb{R}^{4}$ be a complete connected self-shrinker with
$|\vec{H}|^{2} \equiv c$. Then $c\in\{0,1,2\}$. If $c=0$ or $c=1$, then $K\equiv0$. If $c=2$, then $x^{\top}\equiv0$.
\end{proposition}

\begin{proof}

\noindent
Since $x^{\perp}=-\vec{H}$, we have
$|x|^2=|T|^2+c$.
Define
$$u=|x|^{2}-c=|T|^{2}\geq 0.$$
The constancy of $c$ gives
\begin{equation}\label{3.1-1}
\nabla u=2T, \quad |\nabla u|^2=4u.
\end{equation}
Using $\Delta x=\vec{H}$, we also obtain
\begin{equation}\label{3.1-2}
\Delta u=4+2\langle x, \vec{H}\rangle=4-2c.
\end{equation}
On $\Omega=\{u>0\}$, put $r = \sqrt{u}$. Equations \eqref{3.1-1} and \eqref{3.1-2} imply
\begin{equation}\label{3.1-3}
|\nabla r|=1, \quad \Delta r = \frac{1-c}{r}.
\end{equation}
Differentiating $|\nabla r|^{2}=1$ gives $\nabla^{2}r(\nabla r,\cdot)=0$. Thus, in a local orthonormal frame $e_{1}=\nabla r$, $e_{2}\perp e_{1}$, the only possibly nonzero component of $\nabla^{2}r$ is
$$\nabla^{2} r(e_{2}, e_{2})=\frac{1-c}{r}.$$
Consequently,
\begin{equation}\label{3.1-4}
|\nabla^{2} r|^2=\frac{(1-c)^{2}}{r^{2}}.
\end{equation}
The Bochner formula and $\operatorname{Ric}=Kg$ now give
\begin{equation*}
0=\frac{1}{2} \Delta |\nabla r|^{2}
=|\nabla^{2} r|^{2}+\langle \nabla r, \nabla \Delta r\rangle+K|\nabla r|^{2}
=\frac{(1-c)^{2}}{r^{2}}-\frac{1-c}{r^{2}}+K.
\end{equation*}
Hence,
\begin{equation}\label{3.1-5}
K=\frac{c(1-c)}{r^2} \quad \text{on } \Omega.
\end{equation}

Suppose first that $\Omega\neq\varnothing$. Fix $p\in\Omega$ and use Lemma \ref{lemma 2.1} to obtain a curve $\gamma(t)$ converging to $q\in\{u=0\}$ as $t\uparrow r(p)$. Since $K$ is smooth, it is bounded in a sufficiently small neighborhood of $q$. Therefore \eqref{3.1-5} yields
$$c(1-c)=K(\gamma(t))\cdot r(\gamma(t))^{2}\to 0.$$
It follows that $c=0$ or $c=1$. In either case, \eqref{3.1-2} gives $\Delta u>0$, so $\{u=0\}$ has empty interior. Thus $\Omega$ is dense, and \eqref{3.1-5} implies $K\equiv0$ by continuity.

If $\Omega=\varnothing$, then $u\equiv0$. Equation \eqref{3.1-2} gives $c=2$, while $u=|T|^2$ gives $T\equiv0$. Conversely, the first case excludes $c=2$ whenever $\Omega$ is nonempty. This proves the proposition.
\end{proof}

\vskip3mm
\noindent
{\it Proof of Theorem \ref{theorem 1.1}}.

Let $J$ be the standard complex structure on $\mathbb{R}^{4}$. By Proposition \ref{proposition 3.1}, only the three values $c=0,1,2$ need be considered.
If $c=0$ or $c=1$, we know that $K\equiv0$.
Suppose now that $c=2$. Then $x^{\top}=0$, $|x|^{2}=2$, and $\vec{H}=-x$. Since $x$ is normal and the immersion is Lagrangian, $Jx$ is tangent. Define the unit tangent field
$$
E_{1}=\frac{Jx}{\sqrt{2}}.
$$
For every tangent vector $X$, $D_{X}E_{1}=\frac{JX}{\sqrt{2}}$ is normal. Taking tangential and normal components gives
\begin{equation*}
\nabla_{X}E_{1}=0, \quad B(X,E_{1})=\frac{JX}{\sqrt{2}}.
\end{equation*}
In particular, $M^{2}$ has a nonvanishing parallel tangent field, which implies $K\equiv0$.
From the above analysis, it can be seen that
\begin{equation*}
K\equiv0, \quad |B|^{2}=|\vec{H}|^{2}\in\{1, 2, 3\}.
\end{equation*}
Consequently, the conclusion follows from the result of Cheng, Hori and Wei \cite{CHW}.
\begin{flushright}
$\square$
\end{flushright}

\vskip3mm
\noindent
{\it Proof of Theorem \ref{theorem 1.2}}.

Proposition \ref{proposition 3.1} reduces the argument to three cases.

\smallskip\noindent\emph{Case 1: $c=0$.}
Here $K\equiv0$, so the Gauss equation gives $|B|^{2}=0$. The immersion is totally geodesic and $M^{2}$ is isometric to a plane through the origin.

\smallskip\noindent\emph{Case 2: $c=1$.}
The Gauss equation gives $|B|^{2}\equiv1$. On the dense open set $\Omega$, \eqref{3.1-3} implies $\nabla^{2}r=0$. Consequently, $e_{1}=\nabla r$ is parallel, and a local orthogonal unit field $e_{2}$ is also parallel. Since $T=re_{1}$, we have
\begin{equation}\label{3.1-6}
\nabla_{e_{1}}T = e_{1}, \quad \nabla_{e_{2}}T = 0.
\end{equation}
Using the self-shrinker equation in \eqref{2.1-1}, we obtain
\begin{align}\label{3.1-7}
\nabla_{X} T=X-A_{\vec{H}} X, \quad \nabla^{\perp}_{X}\vec{H}=B(X, T).
\end{align}
Combining \eqref{3.1-6} and \eqref{3.1-7} yields
\begin{equation}\label{3.1-8}
A_{\vec{H}}e_{1}=0, \quad A_{\vec{H}}e_{2} = e_{2}, \quad |A_{\vec{H}}|^{2}=1.
\end{equation}
Equation \eqref{2.1-3}, with $|\vec{H}|^{2}=1$, then gives $\nabla^{\perp} H=0$ on $\Omega$, and hence on all of $M$ by continuity.
Then, \eqref{3.1-7} now implies
\begin{equation}\label{3.1-9}
B(e_{1}, e_{1})=0, \quad B(e_{1}, e_{2})=0, \quad B(e_{2}, e_{2})=\vec{H}.
\end{equation}
By use of $|B|^{2}\equiv 1$, \eqref{2.1-5} and \eqref{3.1-9}, we get that
\begin{equation*}
|\nabla^{\perp}B|=0.
\end{equation*}
By density, $\nabla^{\perp}B=0$ on $M$. It follows that $B$ is parallel everywhere, with eigenvalues $0$ and $1$.
\eqref{3.1-9} implies $B_{ij}=\langle B_{ij}, \eta\rangle\eta$ for any $i, j, \eta=\frac{\vec{H}}{|\vec{H}|}$. By the Theorem 1 of Yau (\cite{Yau}), $M^{2}$ lies some
$3$-dimensional linear subspace $\mathbb{R}^{3}$.
In Theorem 4 of \cite{Law}, Lawson showed that
every smooth hypersurface with $\nabla^{\perp}B=0$ splits isometrically as a product of a sphere
and a linear space. That is, $M^{2}$ is isometric to $S^{1}(1)\times\mathbb{R}^{1}$.

\smallskip\noindent\emph{Case 3: $c=2$.}
Proposition \ref{proposition 3.1} gives $T=0$, and hence $|x|^2=2$. Thus $x$ takes values in $S^{3}(\sqrt{2})$. If $\vec{H}^{S}$ denotes its mean curvature vector in that sphere, then
$$
\vec{H}=\vec{H}^{S}-\frac{2}{(\sqrt2)^2}x=\vec{H}^{S}-x.
$$
On the other hand, $\vec{H}=-x$ by the self-shrinker equation, so $\vec{H}^S=0$. That is, $x$ is a minimal immersion in the sphere.
\begin{flushright}
$\square$
\end{flushright}

\noindent{\bf Acknowledgements.}
The first author was partially supported by NSFC grant No.12401060, Natural Science Foundation of Henan Province grant No.262300421232. The second author was partially supported by NSFC Grant No.12171164. The third author was partially supported by NSFC grant No. 12501063, Postdoctoral Fellowship Program of CPSF
grant No. GZC20252036.

\end{document}